\documentclass{article}%
\usepackage{amsmath}%
\usepackage{amsfonts}%
\usepackage{amssymb}%
\usepackage{graphicx}

\begin{document}
\title{On Transforming Functions of a Certain Dot-Product Gradient Operator}
\author{Henrik Stenlund\thanks{The author is grateful to Visilab Signal Technologies for supporting this work.}}
\date{July 21st, 2021}
\maketitle
\begin{abstract}
In this paper it is shown that a function of the constant dot product of the gradient operator acting on an arbitrary function can be transformed to a double three-dimensional integral. The inner one of them is a Fourier transform of the operator function. The result converted to one-dimensional problems is also useful in transforming complex differential expressions.
\footnote{Visilab Report \#2021-07}
\end{abstract}
\subsection{Mathematical Classification}
Mathematics Subject Classification: 34A05, 34A55, 34A30, 34C20, 34L40, 35J05, 35J91
\subsection{Keywords}
Gradient operator, function of gradient, wave equation, three-dimensional differential equations, dot product of gradient
\section{Introduction}
The gradient operator is a bit troublesome to manage in various expressions including differential equations. More difficulties are faced when the gradient is an argument of a function. The case of a function with a vector gradient operator (defined as a power series) is forbiddingly complex as it immediately leads to mixed-rank equations which require completely different methods for progressing. Therefore one is here restricted to scalar-type gradient operators. The exponential operator in the following is defined as a power series and all practical operations with it are based on that. Some cases have already been solved, like the following \cite{Stenlund2021}
\begin{equation}
e^{\alpha\nabla^2}c(\bar{r})=\frac{1}{(4{\pi}\alpha)^{\frac{3}{2}}}\int{d\bar{r'}e^{\frac{-|\bar{r}-\bar{r'}|^2}{4\alpha}}c(\bar{r'})} \label{eqn200}
\end{equation}
\begin{equation}
f({\nabla^2})c(\bar{r})=\int{\frac{d\bar{r'}c(\bar{r'})}{(2\pi)^{3}}}\int{d\bar{k}{e^{i\bar{k}\cdot({\bar{r}}-{\bar{r'}})}f(-k^2)}} \label{eqn2200} 
\end{equation}
and \cite{Stenlund2017}
\begin{equation}
e^{{\beta}\bar{r}\cdot{\nabla}}c(\bar{r})=c(\bar{r}e^{\beta})
\end{equation}
Here the $\alpha$, $\beta$ and $\bar{\beta}$ are complex constant scalars and vectors. It is commonly known \cite{Stenlund2016} that
\begin{equation}
e^{\bar{\beta}\cdot{\nabla}}c(\bar{r})=c(\bar{r}+\bar{\beta}) \label{eqn100} 
\end{equation}
This leads to question, to which form may the following more general expression be converted
\begin{equation}
f(\bar{\beta}\cdot{\nabla})c(\bar{r}) 
\end{equation}
While transforming differential expressions and solving wave equations and other three-dimensional differential equations, the results above and also those obtained in this paper may prove to be useful. The methods shown in \cite{Stenlund2021} are  applied in the following. Any formal proofs are omitted to increase clarity.
\section{The Function of a Gradient Operator Constant Dot Product}
We assume that the function  $c(\bar{r})$ is a Fourier integrable scalar function and may be complex-valued. The function $f()$ is also a Fourier integrable function with some additional requirements in the following. The purpose is to find an expression for the following operator equation
\begin{equation}
f(\bar{\beta}\cdot{\nabla})c(\bar{r}) 
\end{equation}
The three-dimensional Fourier transform of $c(\bar{r})$ is used to get
\begin{equation}
f(\bar{\beta}\cdot{\nabla})\int{\frac{d\bar{k}}{(2\pi)^{\frac{3}{2}}}e^{i\bar{k}\cdot{\bar{r}}}\tilde{c}(\bar{k})} 
\end{equation}
Now it is assumed that the function $f()$ has a MacLaurin series expansion which is a fair assumption for many functions. Using Taylor series will lead to the same end result.
\begin{equation}
f(z)=\sum_{n=0}{\frac{f_n{z}^n}{n!}}
\end{equation}
This expression must be used cautiously since the argument $z$ is a differential operator. The function should not contain dependence of $\bar{r}$ in $f_n$ in such a way that it is compromising the operator due to non-commutativity of these items. Any spatial dependence is preferred to be on the left side of the operator. 

The series expansion should also be converging. In order to determine the convergence of the series, the terms should be studied while the operator is acting on the target function. 
\begin{equation}
f(\bar{\beta}\cdot{\nabla})c(\bar{r})=\sum_{n=0}{\frac{f_n{(\bar{\beta}\cdot{\nabla})}^n}{n!}c(\bar{r})}
\end{equation}
Thus the target function becomes part of the convergence. If the operator function can not be expanded as a converging power series, then the formal results may be false. 

By applying the operator function to the exponential function, one is getting
\begin{equation}
f({\bar{\beta}\cdot{\nabla}}){e^{i\bar{k}\cdot\bar{r}}}=\sum_{n=0}{\frac{f_n\cdot{(\bar{\beta}\cdot{\nabla})}^n}{n!}{e^{i\bar{k}\cdot\bar{r}}}}
\end{equation}
The simple result below is straightforward to prove
\begin{equation}
({\bar{\beta}\cdot{\nabla}}){e^{i\bar{k}\cdot\bar{r}}}=({i\bar{k}\cdot{\bar{\beta}}}){e^{i\bar{k}\cdot\bar{r}}}
\end{equation}
and it yields with the series expansion above
\begin{equation}
f({\bar{\beta}\cdot{\nabla}}){e^{i\bar{k}\cdot\bar{r}}}=f({i\bar{k}\cdot{\bar{\beta}}}){e^{i\bar{k}\cdot\bar{r}}}
\end{equation}
It is assumed that the power series will converge. Thus
\begin{equation}
f({\bar{\beta}\cdot{\nabla}})c(\bar{r})=\int{\frac{d\bar{k}}{(2\pi)^{\frac{3}{2}}}e^{i\bar{k}\cdot{\bar{r}}}f({i\bar{k}\cdot{\bar{\beta}}})\tilde{c}(\bar{k})}  
\end{equation}
By reinserting the original transform one will get
\begin{equation}
f({\bar{\beta}\cdot{\nabla}})c(\bar{r})=\int{d\bar{r'}c(\bar{r'})\int{\frac{d\bar{k}}{(2\pi)^{3}}e^{i\bar{k}\cdot({\bar{r}}-{\bar{r'}})}f({i\bar{k}\cdot{\bar{\beta}}})}} \label{eqn2220} 
\end{equation}
The operator function is converted to a double three-dimensional integral over the function's domain while acting on the target function. The inner integral is a three-dimensional Fourier transform of $f({i\bar{k}\cdot{\bar{\beta}}})$. This result is new. 

The exponential case (\ref{eqn100}) is produced from this equation. That can be observed by substituting the exponential function to the inner integral and recognizing that it is a three-dimensional Dirac delta function. That will push out the $c(\bar{r})$ with an offset $\bar{\beta}$ in the argument. 
\section{Examples of Functions}
In the following are presented some simple examples of applying the method presented. The expressions are flattened to one dimension for simplicity.  No attempt is made towards proving convergence of the operator function. It is \textbf{assumed} that convergence is valid in all respects required since one does not possibly know what the target function $c(x)$ is. 

The result (\ref{eqn2220}) above is flattened to one dimension and as there is no dependence on $y,z$ the inner integrals become Dirac delta functions
\begin{equation}
f({\beta\frac{\partial}{\partial{x}}})c(x)=\int_{-\infty}^\infty{\frac{dx'c(x')}{2\pi}\int_{-\infty}^\infty{dke^{ik(x-x')}f(ik\beta)}} \label{eqn3000} 
\end{equation}
\subsection{A Trivial Example of a Function}
The example at hand is the following operator acting on the function $c()$
\begin{equation}
f({\beta\frac{\partial}{\partial{x}}})=\frac{1}{\beta\frac{\partial}{\partial{x}}} \label{eqn3050} 
\end{equation}
The result would be expected to be some sort of an integral. One can directly substitute the function to the inner integral
\begin{equation}
\int_{-\infty}^\infty{\frac{dx'c(x')}{2\pi}\int_{-\infty}^\infty{\frac{dke^{ik(x-x')}}{ik\beta}}}  
\end{equation}
The inner integral is a tabulated Fourier transform \cite{Gradshteyn2007} yielding
\begin{equation}
\frac{1}{2\beta}\int_{-\infty}^\infty{\frac{dx'c(x')}{2\pi}sgn(x-x')} \label{eqn3090} 
\end{equation}
or
\begin{equation}
\frac{1}{\beta\frac{\partial}{\partial{x}}}c(x)=\frac{1}{2\beta}[-\int_{x}^\infty{dx'c(x')}+\int_{-\infty}^x{dx'c(x')}] \label{eqn3100} 
\end{equation}
By applying the original operator's inverse
\begin{equation}
{\beta\frac{\partial}{\partial{x}}} \label{eqn3350} 
\end{equation}
the result is verified.
\subsection{A Slightly More Complex Example}
By using the result (\ref{eqn3000}) again for the case of
\begin{equation}
f({\beta\frac{\partial}{\partial{x}}})=e^{\beta\frac{\partial}{\partial{x}}}\beta\frac{\partial}{\partial{x}}
\end{equation}
After substitution to the inner integral one will obtain
\begin{equation}
\int_{-\infty}^\infty{\frac{dx'c(x')}{2\pi}\int_{-\infty}^\infty{dke^{ik(x-x'+\beta)}}{ik\beta}} \label{eqn4440}
\end{equation}
The inner integral is processed as follows
\begin{equation}
-\beta\int_{-\infty}^\infty{dke^{ik(x'-x-\beta)}ik}
\end{equation}
where the integration variable has been changed. According to the delta function's definition, one will get it in the form
\begin{equation}
-2\pi\beta\frac{\partial}{\partial{x'}}\delta(x'-x-\beta)
\end{equation}
Placing this back to equation(\ref{eqn4440}) yields
\begin{equation}
-\beta\int_{-\infty}^\infty{{dx'c(x')}\frac{\partial}{\partial{x'}}\delta(x'-x-\beta)} \label{eqn4450}
\end{equation}
Elementary operations offer the following result
\begin{equation}
e^{\beta\frac{\partial}{\partial{x}}}\beta\frac{\partial}{\partial{x}}c(x)=[\beta\frac{\partial}{\partial{x}}c(x)]_{x=x+\beta}
\end{equation}
The result is not surprising and is obvious from what has been shown earlier in this paper and is a combination of two commutative operators.
\subsection{A More Complex Example}
Cases where the function $f({\beta\frac{\partial}{\partial{x}}})$ is a trigonometric function of first order in the argument usually leads to groups of delta functions. These, in turn, create groups of the $c(x)$ at various points. More complex cases come to light when the argument is of higher order. By applying the result (\ref{eqn3000}) for the case of
\begin{equation}
f({\beta\frac{\partial}{\partial{x}}})c(x)=cos(({\beta\frac{\partial}{\partial{x}}})^2)c(x)   \label{eqn6450}
\end{equation}
One can use a tabulated Fourier transform \cite{Gradshteyn2007}
\begin{equation}
\int_{-\infty}^\infty{dk\cdot{e^{ik\rho}cos(bk^2)}}=\sqrt\frac{\pi}{b}cos(\frac{\rho^2}{4b}-\frac{\pi}{4})
\end{equation}
to handle the inner integral, obtaining
\begin{equation}
\int_{-\infty}^\infty{dk\cdot{e^{ik(x-x')}cos(\beta^2{k^2})}}=\frac{\sqrt{\pi}}{\beta}cos[\frac{(x-x')^2}{4\beta^2}-\frac{\pi}{4}]
\end{equation}
The expression ({\ref{eqn6450}) is transformed to the form
\begin{equation}
cos(({\beta\frac{\partial}{\partial{x}}})^2)c(x)=\frac{1}{\sqrt{4\pi}{\beta}}\int_{-\infty}^\infty{dx'\cdot{c(x')}cos[\frac{(x-x')^2}{4\beta^2}-\frac{\pi}{4}]}
\end{equation}
The result is interesting as it shows a very violent behavior of the integrand since the $cos()$ will vary between $+1$ and $-1$ with a progressively growing argument. It is modulating the target function and proves that even very strongly acting differential operators may have a corresponding integral representation.
\section{Conclusions}
This paper presents a result regarding the constant dot product gradient operator as an argument of an arbitrary function. The equation (\ref{eqn2220}) shows how to change a function of a constant dot product gradient operator function to an integral operator with a kernel. To derive this method, Fourier analysis in three dimensions is applied in a straightforward way. The kernel is a Fourier transform of $f({i\bar{k}\cdot{\bar{\beta}}})$. 

The basic requirements for the presented method to be valid are that the functions $c(\bar{r})$ and $f(z)$ are Fourier integrable and that $f({\bar{\beta}\cdot{\nabla}})c(\bar{r})$ is converging when $f(z)$ is expressed as a series.

The new result is applicable in transforming and solving differential expressions for more complex cases than polynomials. The flattened result (\ref{eqn3000}) can be used for transforming one-dimensional differential expressions as long as they do not have interfering $x$ dependence. A few examples are given in terms of it.

\end{document}